# Defining the Mean of a Real-Valued Function on an Arbitrary Metric Space

Kerry M. Soileau

August 9, 2008

### ABSTRACT

We show how a metric space induces a linear functional (a "mean") on real-valued functions with domains in that metric space. This immediately induces a "relative" measure on a collection of subsets of the underlying set.

## 1. INTRODUCTION

In [1], Dembski gave a definition of uniform probability distribution of points in a compact subset of an arbitrary metric space. In this paper we show how this definition can be extended to that of the mean of a function over a compact domain in such a space. Additionally, by regarding the mean of the characteristic function of a set $A$ over a set $B$, we show how a metric space induces the "relative measure" of $A$ over a set $B$. Relative measure is analogous to conditional probability; the relative measure of a set with respect to itself is always unity, and in general the relative measure of $A$ over $B$ may be regarded as the probability that a uniformly randomly chosen point in $B$ also lies in $A$.

## 2. THE MEAN OF A FUNCTION

Let $K$ be compactly-bounded, i.e. a nonempty subset of some compact subset in a metric space $X$ with metric $\delta$. <u>Definition 2.1</u>: For any $\varepsilon > 0$, we say that a set $S \subseteq X$ is an $\varepsilon$-<u>dispersion</u> if and only if for every distinct $x, y \in S$ we have $\delta(x, y) \geq \varepsilon$. We say that $S$ is an $\varepsilon$-<u>lattice</u> of a set $T$ if no set $U$ satisfying $S \subset U \subseteq T$ is an $\varepsilon$-dispersion.

# DEFINING THE MEAN OF A REAL-VALUED FUNCTION ON AN ARBITRARY METRIC SPACE

Let $L_\varepsilon(K)$ be the collection of sets which are $\varepsilon$-lattices of $K$. Since $K$ is compactly-bounded, each member of $L_\varepsilon(K)$ is a finite set. (see Dembski [1]).

Now let $f: K \to \mathbb{R}$ be a function. For each nonempty finite set $S \subseteq K$ we define for $f$ the <u>mean</u> $M_f(S) = \frac{1}{|S|} \sum_{s \in S} f(s)$.

Now we define

$$l_\varepsilon(f) = \inf_{S \in L_\varepsilon(K)} M_f(S)$$

and

$$u_\varepsilon(f) = \sup_{S \in L_\varepsilon(K)} M_f(S).$$

Clearly $l_\varepsilon(f) \leq u_\varepsilon(f)$.

<u>Definition 2.2</u>: We say that $f$ <u>has a mean</u> if and only if $\lim_{\varepsilon \to 0} l_\varepsilon(f) = \lim_{\varepsilon \to 0} u_\varepsilon(f)$. If so we define the common limit $\hat{f}$ as the <u>mean</u> of the function $f$.

We first observe that not every function has a mean. For example, take $f: [0,1] \to \mathbb{R}$ given by

$$f(x) = \begin{cases} 1 & \text{if } x \in [0,1] \cap \mathbb{Q} \\ 0 & \text{if } x \in [0,1] \sim \mathbb{Q} \end{cases}.$$

Then it is easy to see that $l_\varepsilon(f) = 0$ and $u_\varepsilon(f) = 1$ for every $\varepsilon > 0$, thus this function has no mean. On the other hand any constant function has a mean just equal to the constant value of the function (as we would expect intuitively).

# DEFINING THE MEAN OF A REAL-VALUED FUNCTION ON AN ARBITRARY METRIC SPACE

<u>Proposition 2.1</u>: If $f : K \to \mathbb{R}$ has a mean, then so does $-f$, and $\widehat{(-f)} = -\hat{f}$. Here $-f : K \to \mathbb{R}$ is given by $(-f)(x) \equiv -f(x)$.

<u>Proof</u>: Fix $\varepsilon > 0$. Then

$$l_\varepsilon(-f) = \inf_{S \in L_\varepsilon(K)} \frac{1}{|S|} \sum_{s \in S} -f(s) = \inf_{S \in L_\varepsilon(K)} -\frac{1}{|S|} \sum_{s \in S} f(s) = -\sup_{S \in L_\varepsilon(K)} \frac{1}{|S|} \sum_{s \in S} f(s) = -u_\varepsilon(f)$$

and

$$-l_\varepsilon(f) = -\inf_{S \in L_\varepsilon(K)} \frac{1}{|S|} \sum_{s \in S} f(s) = \sup_{S \in L_\varepsilon(K)} -\frac{1}{|S|} \sum_{s \in S} f(s) = \sup_{S \in L_\varepsilon(K)} \frac{1}{|S|} \sum_{s \in S} -f(s) = u_\varepsilon(-f)$$

Hence $l_\varepsilon(-f) = -u_\varepsilon(f)$ and $-l_\varepsilon(f) = u_\varepsilon(-f)$. Taking the limits, we get

$$\lim_{\varepsilon \to 0} l_\varepsilon(-f) = -\lim_{\varepsilon \to 0} u_\varepsilon(f) = -\hat{f} = -\lim_{\varepsilon \to 0} l_\varepsilon(f) = \lim_{\varepsilon \to 0} u_\varepsilon(-f), \text{ thus } -f \text{ has a mean and}$$

$$\widehat{(-f)} = -\hat{f}.$$

<u>Proposition 2.2</u>: If $f \geq \alpha$ has a mean and $\alpha$ is a constant, then $\hat{f} \geq \alpha$.

<u>Proof</u>: Fix $\varepsilon > 0$. Then $l_\varepsilon(f) = \inf_{S \in L_\varepsilon(K)} \frac{1}{|S|} \sum_{s \in S} f(s) \geq \inf_{S \in L_\varepsilon(K)} \frac{1}{|S|} \sum_{s \in S} \alpha = \inf_{S \in L_\varepsilon(K)} \alpha = \alpha$.

Similarly $u_\varepsilon(f) = \sup_{S \in L_\varepsilon(K)} \frac{1}{|S|} \sum_{s \in S} f(s) \geq \sup_{S \in L_\varepsilon(K)} \frac{1}{|S|} \sum_{s \in S} \alpha = \sup_{S \in L_\varepsilon(K)} \alpha = \alpha$, hence

$l_\varepsilon(f) \geq \alpha$ and $u_\varepsilon(f) \geq \alpha$. Taking the limit yields $\lim_{\varepsilon \to 0} l_\varepsilon(f) \geq \alpha$ and $\lim_{\varepsilon \to 0} u_\varepsilon(f) \geq \alpha$, hence $\hat{f} \geq \alpha$.

<u>Corollary 2.1</u>: If $f \leq \alpha$ has a mean and $\alpha$ is a constant, then $\hat{f} \leq \alpha$.

# DEFINING THE MEAN OF A REAL-VALUED FUNCTION ON AN ARBITRARY METRIC SPACE

Proof: Since $f$ has a mean, we have by Proposition 2.1 that $-f$ has a mean and $\widehat{(-f)} = -\hat{f}$. Since $f \leq \alpha$, we have $-f \geq -\alpha$. Thus by Proposition 2.2, we have that $\widehat{(-f)} \geq -\alpha$, hence $-\hat{f} \geq -\alpha$ and the Corollary is proved.

Proposition 2.3: If $f: K \to \mathbb{R}$ has a mean and $\alpha$ is a constant, then $\alpha f$ has a mean and $\widehat{(\alpha f)} = \alpha \hat{f}$. Here $\alpha f : K \to \mathbb{R}$ is given by $(\alpha f)(x) \equiv \alpha f(x)$.

Proof: First we consider the case in which $\alpha \geq 0$. Fix $\varepsilon > 0$. Then

$$l_\varepsilon(\alpha f) = \inf_{S \in L_\varepsilon(K)} \frac{1}{|S|} \sum_{s \in S} \alpha f(s) = \alpha \inf_{S \in L_\varepsilon(K)} \frac{1}{|S|} \sum_{s \in S} f(s) = \alpha l_\varepsilon(f).$$ Similarly

$$u_\varepsilon(\alpha f) = \sup_{S \in L_\varepsilon(K)} \frac{1}{|S|} \sum_{s \in S} \alpha f(s) = \alpha \sup_{S \in L_\varepsilon(K)} \frac{1}{|S|} \sum_{s \in S} f(s) = \alpha u_\varepsilon(f),$$ hence $l_\varepsilon(\alpha f) = \alpha l_\varepsilon(f)$

and $\alpha u_\varepsilon(f) = u_\varepsilon(\alpha f)$. Taking limits yields

$$\lim_{\varepsilon \to 0} l_\varepsilon(\alpha f) = \alpha \lim_{\varepsilon \to 0} l_\varepsilon(f) = \alpha \hat{f} = \alpha \lim_{\varepsilon \to 0} u_\varepsilon(f) = \lim_{\varepsilon \to 0} u_\varepsilon(\alpha f),$$ so $\alpha f$ has a mean and

$\widehat{(\alpha f)} = \alpha \hat{f}$.

Now if $\alpha < 0$, we have $\widehat{(\alpha f)} = \widehat{((-\alpha)(-f))} = (-\alpha)\widehat{(-f)} = (-\alpha)(-\hat{f}) = \alpha \hat{f}$, and the

Proposition is proved.

Proposition 2.4: If $f, g : K \to \mathbb{R}$ have means, then $f + g$ has a mean and

$\widehat{(f+g)} = \hat{f} + \hat{g}$. Here $(f+g) : K \to \mathbb{R}$ is given by $(f+g)(x) \equiv f(x) + g(x)$.

Proof: Fix $\varepsilon > 0$. Then

# DEFINING THE MEAN OF A REAL-VALUED FUNCTION ON AN ARBITRARY METRIC SPACE

$$l_\varepsilon(f+g) = \inf_{S \in L_\varepsilon(K)} \frac{1}{|S|} \sum_{s \in S} (f+g)(s) = \inf_{S \in L_\varepsilon(K)} \frac{1}{|S|} \sum_{s \in S} (f(s)+g(s))$$

$$\geq \inf_{S \in L_\varepsilon(K)} \frac{1}{|S|} \sum_{s \in S} f(s) + \inf_{S \in L_\varepsilon(K)} \frac{1}{|S|} \sum_{s \in S} g(s) = l_\varepsilon(f) + l_\varepsilon(g)$$

Similarly

$$u_\varepsilon(f+g) = \sup_{S \in L_\varepsilon(K)} \frac{1}{|S|} \sum_{s \in S} (f+g)(s) = \sup_{S \in L_\varepsilon(K)} \frac{1}{|S|} \sum_{s \in S} (f(s)+g(s))$$

$$\leq \sup_{S \in L_\varepsilon(K)} \frac{1}{|S|} \sum_{s \in S} f(s) + \sup_{S \in L_\varepsilon(K)} \frac{1}{|S|} \sum_{s \in S} g(s) = u_\varepsilon(f) + u_\varepsilon(g)$$

so $l_\varepsilon(f) + l_\varepsilon(g) \leq l_\varepsilon(f+g) \leq u_\varepsilon(f+g) \leq u_\varepsilon(f) + u_\varepsilon(g)$. Taking limits yields

$\hat{f} + \hat{g} \leq \lim_{\varepsilon \to 0} l_\varepsilon(f+g) \leq \lim_{\varepsilon \to 0} u_\varepsilon(f+g) \leq \hat{f} + \hat{g}$, hence $f + g$ has a mean and

$\widehat{(f+g)} = \hat{f} + \hat{g}$.

<u>Proposition 2.5</u>: If $f, g : K \to \mathbb{R}$ have means and $\alpha, \beta$ are constants, then

$\widehat{(\alpha f + \beta g)} = \alpha \hat{f} + \beta \hat{g}$.

<u>Proof</u>: $\alpha f$ and $\beta g$ have means, hence $\widehat{(\alpha f + \beta g)} = \widehat{(\alpha f)} + \widehat{(\beta g)} = \alpha \hat{f} + \beta \hat{g}$.

<u>Proposition 2.6</u>: If $f, g : K \to \mathbb{R}$ have means and $f \geq g$, then $\hat{f} \geq \hat{g}$.

<u>Proof</u>: Fix $\varepsilon > 0$. Then $l_\varepsilon(f) = \inf_{S \in L_\varepsilon(K)} \frac{1}{|S|} \sum_{s \in S} f(s) \geq \inf_{S \in L_\varepsilon(K)} \frac{1}{|S|} \sum_{s \in S} g(s) = l_\varepsilon(g)$, hence

$l_\varepsilon(f) \geq l_\varepsilon(g)$. Taking the limit yields $\lim_{\varepsilon \to 0} l_\varepsilon(f) \geq \lim_{\varepsilon \to 0} l_\varepsilon(g)$, hence $\hat{f} \geq \hat{g}$.

<u>Definition 2.3:</u> We say that $K$ is <u>regular</u> if $\liminf_{\varepsilon \to 0} \{|S|; S \in L_\varepsilon(K)\} = \infty$.

Note: Not all $K$ are regular; in particular if $K$ is finite, $K$ is not regular. On the other hand, the closed interval $[0,1]$ is regular.

# DEFINING THE MEAN OF A REAL-VALUED FUNCTION ON AN ARBITRARY METRIC SPACE

Lemma 2.1: If $f : K \to \mathbb{R}$, $K$ is regular and $f = 0$ except at finitely many points, then $f$ has a mean and $\hat{f} = 0$.

Proof: We first prove the case in which $f = 0$ except at a single point. We can assume without loss of generality that $f = 1$ at this unique point. Fix $\tau > 0$. Then for any sufficiently small $\varepsilon > 0$, we can find $\overline{S} \in L_\varepsilon(K)$ such that $\frac{1}{\tau} \leq |\overline{S}|$, hence

$$l_\varepsilon(f) = \inf_{S \in L_\varepsilon(K)} \frac{1}{|S|} \sum_{s \in S} f(s) \leq \frac{1}{|\overline{S}|} \sum_{s \in \overline{S}} f(s) \leq \frac{1}{|\overline{S}|} \leq \tau.$$ Thus $\lim_{\varepsilon \to 0} l_\varepsilon(f) = 0$. Similarly, for any sufficiently small $\varepsilon > 0$, for every $S \in L_\varepsilon(K)$ we have $\frac{1}{\tau} \leq |S|$, hence

$$u_\varepsilon(f) = \sup_{S \in L_\varepsilon(K)} \frac{1}{|S|} \sum_{s \in S} f(s) \leq \sup_{S \in L_\varepsilon(K)} \frac{1}{|S|} \leq \sup_{S \in L_\varepsilon(K)} \tau = \tau.$$ Thus $\lim_{\varepsilon \to 0} u_\varepsilon(f) = 0$. Hence $f$ has a mean and $\hat{f} = 0$.

Next, we prove the full Lemma inductively. We have just shown that the Lemma is true for $f = 0$ except at a single point. Now assume it is true for $f = 0$ except at $k$ points. Let $g : K \to \mathbb{R}$ be a function such that $g = 0$ except at $k+1$ points. We may express $g$ as the sum $g_1 + g_k$, where $g_1 = 0$ except at a single point, and $g_k = 0$ except at $k$ points. Then $g_1$ has a mean and $\hat{g}_1 = 0$, and $g_k$ has a mean and $\hat{g}_k = 0$. Hence $g$ has a mean and $\hat{g} = \hat{g}_1 + \hat{g}_k = 0$, and the Lemma is proved.

Proposition 2.7: If $f, g : K \to \mathbb{R}$, $f$ has a mean, $K$ is regular and $g = f$ except at finitely many points, then $g$ has a mean and $\hat{g} = \hat{f}$.

# DEFINING THE MEAN OF A REAL-VALUED FUNCTION ON AN ARBITRARY METRIC SPACE

<u>Proof</u>: Write $g = (g-f) + f$. By Proposition 2.5, $g-f$ has a mean. By Lemma 2.1, $g-f$ has a mean of $0$, and $f$ has a mean by assumption so $g = (g-f) + f$ has a mean and $\hat{g} = \widehat{(g-f)} + \hat{f} = 0 + \hat{f} = \hat{f}$, and the Proposition is proved.

<u>Theorem 2.1</u>: If each $f_n : K \to \mathbb{R}$ has a mean, $f_n \to f$ uniformly on $K$, and each $|\hat{f}_n| < \infty$, then $f$ has a mean and $\hat{f} = \lim_{n \to \infty} \hat{f}_n$.

<u>Proof</u>: Fix $\tau > 0$. Since $f_n \to f$ uniformly on $K$, we can find an $n \geq 1$ such that

$$|f_n - f| \leq \frac{\tau}{3} \text{ on } K, \text{ hence } f_n - \frac{\tau}{3} \leq f \leq f_n + \frac{\tau}{3} \text{ there. Next, since } |\hat{f}_n| < \infty, \text{ it follows that}$$

$\lim_{\varepsilon \to 0} |l_\varepsilon(f_n) - u_\varepsilon(f_n)| = 0$, hence we can find $\bar{\varepsilon} > 0$ such that $|l_\varepsilon(f_n) - u_\varepsilon(f_n)| \leq \frac{\tau}{3}$ for all $\varepsilon < \bar{\varepsilon}$. Next, since $f_n - \frac{\tau}{3} \leq f \leq f_n + \frac{\tau}{3}$, we have

$$l_\varepsilon(f_n) - \frac{\tau}{3} = l_\varepsilon\left(f_n - \frac{\tau}{3}\right) \leq l_\varepsilon(f) \leq l_\varepsilon\left(f_n + \frac{\tau}{3}\right) = l_\varepsilon(f_n) + \frac{\tau}{3} \text{ and}$$

$$u_\varepsilon(f_n) - \frac{\tau}{3} = u_\varepsilon\left(f_n - \frac{\tau}{3}\right) \leq u_\varepsilon(f) \leq u_\varepsilon\left(f_n + \frac{\tau}{3}\right) = u_\varepsilon(f_n) + \frac{\tau}{3} \text{ for } \varepsilon < \bar{\varepsilon}. \text{ Equivalently, we}$$

have $|l_\varepsilon(f_n) - l_\varepsilon(f)| \leq \frac{\tau}{3}$ and $|u_\varepsilon(f_n) - u_\varepsilon(f)| \leq \frac{\tau}{3}$ for $\varepsilon < \bar{\varepsilon}$.

Finally

$$|u_\varepsilon(f) - l_\varepsilon(f)| \leq |u_\varepsilon(f) - u_\varepsilon(f_n)| + |u_\varepsilon(f_n) - l_\varepsilon(f_n)| + |l_\varepsilon(f_n) - l_\varepsilon(f)| \leq \frac{\tau}{3} + \frac{\tau}{3} + \frac{\tau}{3} = \tau \text{ for}$$

for $\varepsilon < \bar{\varepsilon}$. Since $\tau$ was an arbitrary positive number, we have shown that

$\lim_{\varepsilon \to 0} |u_\varepsilon(f) - l_\varepsilon(f)| = 0$, $\lim_{\varepsilon \to 0} |l_\varepsilon(f_n) - l_\varepsilon(f)| = 0$ and $\lim_{\varepsilon \to 0} |u_\varepsilon(f_n) - u_\varepsilon(f)| = 0$.

# DEFINING THE MEAN OF A REAL-VALUED FUNCTION ON AN ARBITRARY METRIC SPACE

Now since $\left|\hat{f}_n\right| < \infty$, we have $\left|\lim_{\varepsilon \to 0} l_\varepsilon(f_n)\right| = \left|\hat{f}_n\right| = \left|\lim_{\varepsilon \to 0} u_\varepsilon(f_n)\right| < \infty$, hence for some $M$ we have $\left|l_\varepsilon(f_n)\right| \leq M$ and $\left|u_\varepsilon(f_n)\right| \leq M$ for $\varepsilon < \overline{\varepsilon}$. Thus

$$\left|l_\varepsilon(f)\right| \leq \left|l_\varepsilon(f) - l_\varepsilon(f_n)\right| + \left|l_\varepsilon(f_n)\right| \leq \left|l_\varepsilon(f) - l_\varepsilon(f_n)\right| + M \text{ and}$$

$$\left|u_\varepsilon(f)\right| \leq \left|u_\varepsilon(f) - u_\varepsilon(f_n)\right| + \left|u_\varepsilon(f_n)\right| \leq \left|u_\varepsilon(f) - u_\varepsilon(f_n)\right| + M \text{ . Since } \lim_{\varepsilon \to 0}\left|l_\varepsilon(f_n) - l_\varepsilon(f)\right| = 0$$

and $\lim_{\varepsilon \to 0}\left|u_\varepsilon(f_n) - u_\varepsilon(f)\right| = 0$, $\left|l_\varepsilon(f)\right|$ and $\left|u_\varepsilon(f)\right|$ are clearly bounded for $\varepsilon < \overline{\varepsilon}$. Since

$\lim_{\varepsilon \to 0}\left|u_\varepsilon(f) - l_\varepsilon(f)\right| = 0$, we may conclude that $\left|\hat{f}\right| = \lim_{\varepsilon \to 0}\left|l_\varepsilon(f)\right| = \lim_{\varepsilon \to 0}\left|u_\varepsilon(f)\right| < \infty$.

To complete the proof, note that for any $\tau > 0$ we can find $N$ such that

$f_n - \tau \leq f \leq f_n + \tau$ for every $n \geq N$, hence for such $n$, and for any $\varepsilon > 0$, we have

$u_\varepsilon(f_n) - \tau \leq u_\varepsilon(f) \leq u_\varepsilon(f_n) + \tau$. Taking limits as $\varepsilon \to 0$, we get $\hat{f}_n - \tau \leq \hat{f} \leq \hat{f}_n + \tau$, i.e.

$\left|\hat{f}_n - \hat{f}\right| \leq \tau$. Thus $\hat{f} = \lim_{n \to \infty} \hat{f}_n$.

## 3. RELATIVE MEASURABILITY

<u>Definition 3.1</u>: We say that a set $S \subseteq K$ is <u>measurable relative to</u> $K$ if the function

$\chi_S : K \to \mathbb{R}$ given by

$$\chi_S(x) = \begin{cases} 1 & \text{if } x \in S \\ 0 & \text{if } x \in K \sim S \end{cases}$$

has a mean.

<u>Proposition 3.1</u>: If $A$ is measurable relative to $K$ and $\chi_A : K \to \mathbb{R}$ and $\chi_{K \sim A} : K \to \mathbb{R}$ are

given by

$$\chi_A(x) = \begin{cases} 1 & \text{if } x \in A \\ 0 & \text{if } x \in K \sim A \end{cases}$$

# DEFINING THE MEAN OF A REAL-VALUED FUNCTION ON AN ARBITRARY METRIC SPACE

and

$$\chi_{K\sim A}(x) = \begin{cases} 0 & \text{if } x \in A \\ 1 & \text{if } x \in K \sim A \end{cases},$$

then $\widehat{\chi_A} + \widehat{\chi_{K\sim A}} = 1$.

Proof: Fix $\varepsilon > 0$. Then

$$l_\varepsilon(\chi_A) + l_\varepsilon(\chi_{K\sim A}) = \inf_{S \in L_\varepsilon(K)} \frac{1}{|S|} \sum_{s \in S} \chi_A(s) + \inf_{S \in L_\varepsilon(K)} \frac{1}{|S|} \sum_{s \in S} \chi_{K\sim A}(s)$$
$$\leq \inf_{S \in L_\varepsilon(K)} \frac{1}{|S|} \sum_{s \in S} (\chi_A(s) + \chi_{K\sim A}(s)) = \inf_{S \in L_\varepsilon(K)} \frac{1}{|S|} \sum_{s \in S} 1 = \inf_{S \in L_\varepsilon(K)} 1 = 1.$$

Similarly

$$u_\varepsilon(\chi_A) + u_\varepsilon(\chi_{K\sim A}) = \sup_{S \in L_\varepsilon(K)} \frac{1}{|S|} \sum_{s \in S} \chi_A(s) + \sup_{S \in L_\varepsilon(K)} \frac{1}{|S|} \sum_{s \in S} \chi_{K\sim A}(s)$$
$$\geq \sup_{S \in L_\varepsilon(K)} \frac{1}{|S|} \sum_{s \in S} (\chi_A(s) + \chi_{K\sim A}(s)) = \sup_{S \in L_\varepsilon(K)} \frac{1}{|S|} \sum_{s \in S} 1 = \sup_{S \in L_\varepsilon(K)} 1 = 1,$$

hence

$l_\varepsilon(\chi_A) + l_\varepsilon(\chi_{K\sim A}) \leq 1$ and $u_\varepsilon(\chi_A) + u_\varepsilon(\chi_{K\sim A}) \geq 1$. Taking the limits yields $\widehat{\chi_A} + \widehat{\chi_{K\sim A}} \leq 1$ and $\widehat{\chi_A} + \widehat{\chi_{K\sim A}} \geq 1$, hence $\widehat{\chi_A} + \widehat{\chi_{K\sim A}} = 1$.

Definition 3.2: If $A \subseteq B$ and $A$ is measurable relative to $B$, then we define the <u>measure of $A$ relative to $B$</u> as $\widehat{\chi_{A|B}}$, where

$$\chi_{A|B}(x) = \begin{cases} 1 & \text{if } x \in A \\ 0 & \text{if } x \in B \sim A \end{cases}.$$

We abbreviate $\widehat{\chi_{A|B}}$ as $(A|B)$.

# DEFINING THE MEAN OF A REAL-VALUED FUNCTION ON AN ARBITRARY METRIC SPACE

Definition 3.3: Let $A \subseteq B$. We say that $A$ <u>has a thin boundary with</u> $B$ if $A$ is measurable relative to $B$, and $\lim_{\varepsilon \to 0} \inf_{S \in L_\varepsilon(K)} \frac{|A \cap S|}{|B \cap S|} = (A|B) = \lim_{\varepsilon \to 0} \sup_{S \in L_\varepsilon(K)} \frac{|A \cap S|}{|B \cap S|}$ for any compactly-bounded superset $K \supseteq B$, where $(A|B)$ is a real number with value independent of the choice of $K$. Note that any set has a thin boundary with itself.

Proposition 3.2: If $A$ has a thin boundary with $B$ and $B$ has a thin boundary with $C$, and $C$ is compactly-bounded, then $A$ has a thin boundary with $C$ and

$$(A|C) = (A|B)(B|C).$$

Proof: Fix $\alpha > 0$. Since $\lim_{\tau \to 0} \max \begin{pmatrix} |((A|B) - \tau)((B|C) - \tau) - (A|B)(B|C)|, \\ |((A|B) + \tau)((B|C) + \tau) - (A|B)(B|C)| \end{pmatrix} = 0$, we can find $\tau$ satisfying

$$\max \begin{pmatrix} |((A|B) - \tau)((B|C) - \tau) - (A|B)(B|C)|, \\ |((A|B) + \tau)((B|C) + \tau) - (A|B)(B|C)| \end{pmatrix} < \alpha.$$

Fix a compactly-bounded superset $K \supseteq C$. Since $B$ has a thin boundary with $C$, we have $\lim_{\varepsilon \to 0} \inf_{S \in L_\varepsilon(K)} \frac{|B \cap S|}{|C \cap S|} = (B|C) = \lim_{\varepsilon \to 0} \sup_{S \in L_\varepsilon(K)} \frac{|B \cap S|}{|C \cap S|}$, hence we can find $\bar{\varepsilon} > 0$ such that

$((B|C) - \tau)|C \cap S| < |B \cap S| < ((B|C) + \tau)|S|$ for every positive $\varepsilon < \bar{\varepsilon}$ and $S \in L_\varepsilon(K)$.

Similarly, since $A$ has a thin boundary with $B$, we have

$\lim_{\varepsilon \to 0} \inf_{S \in L_\varepsilon(K)} \frac{|A \cap S|}{|B \cap S|} = (A|B) = \lim_{\varepsilon \to 0} \sup_{S \in L_\varepsilon(K)} \frac{|A \cap S|}{|B \cap S|}$, hence we can find $\bar{\varepsilon} > 0$ such that

$((A|B) - \tau)|B \cap S| < |A \cap S| < ((A|B) + \tau)|B \cap S|$ for every positive $\varepsilon < \bar{\varepsilon}$ and

$S \in L_\varepsilon(K)$. Fix $\bar{S} \in L_\varepsilon(K)$. Then $((B|C) - \tau)|C \cap \bar{S}| < |B \cap \bar{S}| < ((B|C) + \tau)|C \cap \bar{S}|$.

# DEFINING THE MEAN OF A REAL-VALUED FUNCTION ON AN ARBITRARY METRIC SPACE

We also have $((A|B)-\tau)|B\cap \bar{S}| < |A\cap \bar{S}| < ((A|B)+\tau)|B\cap \bar{S}|$. These inequalities imply

$$((A|B)-\tau) < \frac{|A\cap \bar{S}|}{|B\cap \bar{S}|} < ((A|B)+\tau) \text{ and } ((B|C)-\tau) < \frac{|B\cap \bar{S}|}{|C\cap \bar{S}|} < ((B|C)+\tau),$$

hence we have $((A|B)-\tau)((B|C)-\tau) < \frac{|A\cap \bar{S}|}{|B\cap \bar{S}|}\frac{|B\cap \bar{S}|}{|C\cap \bar{S}|} < ((A|B)+\tau)((B|C)+\tau)$ or equivalently

$$((A|B)-\tau)((B|C)-\tau)-(A|B)(B|C) < \frac{|A\cap \bar{S}|}{|C\cap \bar{S}|}-(A|B)(B|C)$$
$$< ((A|B)+\tau)((B|C)+\tau)-(A|B)(B|C).$$

By choice of $\tau$ we may conclude that $\left|\frac{|A\cap \bar{S}|}{|C\cap \bar{S}|}-(A|B)(B|C)\right| < \alpha$. Since $\bar{S}$ was an arbitrary member of $L_\varepsilon(K)$, we have that $-\alpha < \frac{|A\cap S|}{|C\cap S|}-(A|B)(B|C) < \alpha$ holds for every $S \in L_\varepsilon(K)$. Hence $-\alpha < \sup_{S\in L_\varepsilon(K)} \frac{|A\cap S|}{|C\cap S|}-(A|B)(B|C) \leq \alpha$ and

$-\alpha \leq \inf_{S\in L_\varepsilon(K)} \frac{|A\cap S|}{|C\cap S|}-(A|B)(B|C) < \alpha$. Since $\alpha$ was an arbitrary positive number, it follows that $\lim_{\varepsilon \to 0}\inf_{S\in L_\varepsilon(K)} \frac{|A\cap S|}{|C\cap S|} = (A|B)(B|C)$ and $\sup_{S\in L_\varepsilon(K)} \frac{|A\cap S|}{|C\cap S|} = (A|B)(B|C)$.

Hence $A$ has a thin boundary with $C$ and $(A|C) = (A|B)(B|C)$.

Next note that relative measure has an interesting subadditivity property.

# DEFINING THE MEAN OF A REAL-VALUED FUNCTION ON AN ARBITRARY METRIC SPACE

<u>Proposition 3.3</u>: If each member of a disjoint sequence of nonempty sets $\{A_n\}$ has a thin boundary with $K$, and $\lim_{\varepsilon \to \infty} \sum_{n=1}^{\infty} l_\varepsilon(\chi_{A_n}) = \sum_{n=1}^{\infty} \hat{\chi}_{A_n} = \lim_{\varepsilon \to \infty} \sum_{n=1}^{\infty} u_\varepsilon(\chi_{A_n})$, then $\bigcup_{n=1}^{\infty} A_n$ has a thin boundary with $K$ and $\left( \left( \bigcup_{n=1}^{\infty} A_n \right) | K \right) = \sum_{n=1}^{\infty} (A_n | K)$.

<u>Proof</u>: Note that

$$\sum_{n=1}^{\infty} l_\varepsilon(\chi_{A_n}) = \sum_{n=1}^{\infty} \inf_{S \in L_\varepsilon(K)} \frac{|A_n \cap S|}{|S|} \leq \inf_{S \in L_\varepsilon(K)} \sum_{n=1}^{\infty} \frac{|A_n \cap S|}{|S|} = \inf_{S \in L_\varepsilon(K)} \frac{\left| \bigcup_{n=1}^{\infty} A_n \cap S \right|}{|S|}$$

$$\leq \sup_{S \in L_\varepsilon(K)} \frac{\left| \bigcup_{n=1}^{\infty} A_n \cap S \right|}{|S|} = \sup_{S \in L_\varepsilon(K)} \sum_{n=1}^{\infty} \frac{|A_n \cap S|}{|S|} \leq \sum_{n=1}^{\infty} \sup_{S \in L_\varepsilon(K)} \frac{|A_n \cap S|}{|S|} = \sum_{n=1}^{\infty} u_\varepsilon(\chi_{A_n})$$

And upon taking limits as $\varepsilon \to 0$, we get

$$\sum_{n=1}^{\infty} \hat{\chi}_{A_n} \leq \lim_{\varepsilon \to 0} \inf_{S \in L_\varepsilon(K)} \frac{\left| \bigcup_{n=1}^{\infty} A_n \cap S \right|}{|S|} \leq \lim_{\varepsilon \to 0} \inf_{S \in L_\varepsilon(K)} \sup_{S \in L_\varepsilon(K)} \frac{\left| \bigcup_{n=1}^{\infty} A_n \cap S \right|}{|S|} \leq \sum_{n=1}^{\infty} \hat{\chi}_{A_n}$$, thus $\bigcup_{n=1}^{\infty} A_n$ has a thin boundary with $K$ and $\left( \bigcup_{n=1}^{\infty} A_n | K \right) = \sum_{n=1}^{\infty} (A_n | K)$.

International Space Station Program Office, Avionics and Software Office, Mail Code OD, NASA Johnson Space Center, Houston, TX 77058
E-mail address: ksoileau@yahoo.com